\documentclass[12pt, twoside]{amsart}

\swapnumbers\newtheorem{theorem}{Theorem}[section]

\parskip=\medskipamount

\theoremstyle{plain}
\newtheorem{lemma}[theorem]{Lemma}
\newtheorem{prop}[theorem]{Proposition}
\newtheorem{cory}[theorem]{Corollary}
\newtheorem{propdef}[theorem]{Proposition--Definition}

\theoremstyle{definition}
\newtheorem{df}[theorem]{Definition}
\newtheorem{remark}[theorem]{Remark}

\newtheorem{remque}[theorem]{Remarks and Questions}
\newtheorem{example}[theorem]{Example}

\usepackage{amssymb,latexsym, amsmath, amscd, array}

\numberwithin{equation}{section}
\def\CC{\mathbb C}
\def\QQ{\mathbb Q}
\def\QQ{\mathbb Q}
\def\RR{\mathbb R}

\def\ZZ{\mathbb Z}
\def\eps{\varepsilon}

\def\Im{\operatorname {Im}}
\def\Ker{\operatorname {Ker}}
\def\dim{\operatorname {dim}}

\def\rk{\operatorname {rank}}

\def\Aut{\operatorname {Aut}}
\def\hqed{\hfill\hfill$\square$}

\def\sympl{\operatorname {Symp}}
\def\ev{\operatorname {ev}}
\def\ham{\operatorname {Ham}}
\def\Int{\operatorname {Int}}

\newcommand{\ba}{\begin{array}}
\newcommand{\ea}{\end{array}}

\def\AAA{{\mathcal A}}
\def\BBB{{\mathcal B}}

\def\pa{\partial}
\def\la{\langle}
\def\ra{\rangle}

\def\ga{\alpha}
\def\gb{\beta}

\def\gf{\varphi}

\def\gl{\lambda}

\def\gs{\sigma}

\def\CC{\mathbb C}

\def\QQ{\mathbb Q} 
\def\RR{\mathbb R}
\def\ZZ{\mathbb Z}

\def\m{\medskip}
\def\p{{\it Proof.\ }}
\def\wt #1{\widetilde #1}

\begin{document}

\title[Symplectically aspherical manifolds]{Symplectically aspherical 
manifolds}

\author{R. IB\'A\~NEZ}
\author{J. K\c EDRA}
\author{YU. RUDYAK}
\author{A. TRALLE}
\address{R. Ib\'a\~nez, Departamento de Matem\'aticas, Facultad de
Ciencias, Universidad del Pais Vasco, Apdo. 644, 48080 Bilbao, Spain}
\email {mtpibtor@lg.ehu.es}
\address{J. K\c edra, Institute of Mathematics, University of Szczecin,
Wielkopolska 15, 70451 Szczecin, Poland}
\email{kedra@univ.szczecin.pl}
\address {Yu. Rudyak, Department of Mathematics, University of Florida, 
358 Little Hall, Gainesville, FL 32601, USA}
\email{rudyak@math.ufl.edu} 
\address{A. Tralle, Department of Mathematics, University of Warmia and 
Mazury,
10561 Olsztyn, Poland}
\email{tralle@matman.uwm.edu.pl}

\begin{abstract} 
The main subjects of the paper is studying the fundamental groups 
of closed symplectically aspherical manifolds. Motivated by some results of 
Gompf, we  introduce two classes of fundamental groups $\pi_1(M)$ of 
symplectically  aspherical manifolds $M$ with $\pi_2(M)=0$ and $\pi_2(M)\neq 0$. 
Relations between these classes are discussed. We show that several important 
classes of groups can be realized in both classes, while some of groups can be 
realized in the first class but not in the second one. Also, we notice that 
there are some interesting dimensional phenomena in the realization problem.

The above results are framed by a general research of symplectically aspherical 
manifolds. For example, we find some conditions which imply that the 
Gompf sum of symplectically aspherical manifolds is symplectically 
aspherical, or that a total space of a bundle is symplectically aspherical, etc.
\end{abstract}

\subjclass {57R15, 53D05, 14F35} 

\maketitle

\section*{Introduction}
 We say that the symplectic form $\omega$ on a smooth manifold $M$ is {\it 
symplectically aspherical} if
$$
\int_{S^2}f^*\omega=0,
$$ 
for every map $f: S^2\to M$. In cohomological terms, it means that
$$
\langle [\omega], h(a)\rangle=0
$$
for all $a\in\pi_2(M)$, where $h:\pi_2(M)\to H_2(M)$ is the Hurewicz
homomorphism. Frequently one writes the last equality as
$[\omega]|_{\pi_2(M)}=0$. 

\m By the definition, a symplectically aspherical manifold is a connected 
smooth 
manifold which admits a  symplectically aspherical symplectic form. The 
importance of symplectically aspherical  manifolds in symplectic geometry and 
topology is well-known, see e.g. \cite{F, H, LO, R2,
RO, RT}. 
 
\m Clearly, every manifold $M$ with $\pi_2(M)=0$ is symplectically
aspherical whenever it admits a symplectic structure. On the other hand there 
are some reasons to know whether there are symplectically aspherical manifolds 
with $\pi_2(M)\ne 0$, see e.g. \cite{G2}. Examples of such manifolds were 
given 
in \cite{G2} as some 4-dimensional closed manifolds obtained as branched 
coverings.

\m The starting point of this paper was searching for fundamental groups of
symplectically aspherical manifolds. It is well known that every finitely 
presented group can be realized as the fundamental group of a closed 
symplectic 
manifold, \cite{G1}. However, not every such group can be realized as the 
fundamental group of a closed symplectically aspherical manifold. For example, 
the trivial (or, more generally, any finite) group cannot. 

\m
 According to what we have said above, it is interesting to compare the 
fundamental groups of symplectically aspherical manifolds having $\pi_2(M)=0$ 
with these ones having $\pi_2(M)\ne 0$. Definitely, the first class of groups 
is not contained in the second one (for example, the group $\ZZ\oplus \ZZ$ does 
not belong to the second class, cf. Proposition ~\ref{dim} ), but we do not know 
whether the second class is contained in the first one. We also found some 
interesting dimensional phenomena: for instance, the group $\ZZ^{10}$ can be 
realized as the fundamental group of a closed symplectically aspherical 
manifold 
$M^{2n}$ with $2\le n \le 5$ but cannot for $n>5$.

\m However, later we realized that we need more knowledge on symplectically 
aspherical manifolds in themselves. For example, in order to use the van 
Kampen 
Theorem we must know that certain push-outs are symplectically aspherical. 
Because of this, we prove some results on contructions which respect 
symplectic 
asphericity. For example, we found some condition which force the Gompf 
symplectic sum to be symplectically aspherical.

\m So, the paper is divided into 2 parts: In the first one we discuss symplectic 
asphericity of manifolds, in the second one we discuss the fundamental groups 
of symplectically aspherical manifolds. See the table of contents for details.

\m  
We identify de Rham cohomology of a manifold $M$ with $H^*(M;\RR)$. 

Following \cite{LO}, we call a closed connected manifolod $M^{2n}$ {\it 
cohomologically symplectic} or, brielfy, {\it c-symplectic} if there exists a 
class $a\in H^2(M;\RR)$ with $a^n\ne 0$.

{\bf Acknowledgement:}  The paper was supported by MCyT, project BFM 2002-00788, 
Spain. The first author was partially supported by the project UPV/EHU
00127.310-EA-7781/2000, Spain. The second author is a member of European 
Differential Geometry Endeavour (EDGE), Research Training Network
HPRN-CT-2000-00101, supported by The European Human Potential Programme. The 
third author was partially supported by Max-Planck Institute of Mathematics, 
Bonn, Germany, and by the 
free-term research money from the University
of Florida,  Gainesville, USA. The fourth author 
acknowledges the support  of the Polish Committee for the Scientific Research 
(KBN).
 
We are grateful to Vladimir Chernov and Dusa McDuff for useful discussions and 
valuable advices.

\tableofcontents

\section{Preliminaries}
First,  we fix the following fact.

\begin{prop}
\label{induced}
Let $\omega$ be a symplectically aspherical form on a manifold $N$, and let 
$g: 
M \to N$ be a map such that $g^*\omega$ is a symplectic form on $M$. Then $M$ 
is 
symplectically aspherical. In particular, a covering manifold over 
symplectically aspherical manifold is symplectically 
aspherical.
\hqed
\end{prop}

\m We need the following homotopic characterization of symplectically
aspherical manifolds. Given a group $\pi$, recall that the Eilenberg--Mac Lane 
space $K(\pi,1)$ is defined to be a connected $CW$-space with the fundamental 
group $\pi$ and such that $\pi_i(K(\pi,1)=0$ for $i>1$. It follows from the 
elementary obstruction theory that the homotopy type of such space is completely 
determined by $\pi$ 

\begin{prop}  
\label{aspher}
Let $(M,\omega)$ be a smooth manifold manifold, let $\omega$ be a symplectic 
form on $M$, and let $K$ denote the Eilenberg--Mac Lane space $K(\pi_1(M),1)$. 
The following three
conditions are equivalent:
\par {\rm (i)} the form $\omega$ is symplectically aspherical; 

\par {\rm (ii)} there exists a map $f: M\to K$ which induces isomorphism on 
the
fundamental groups and such that 
$$ 
[\omega]\in\Im\,\{f^*: H^2(K,\RR)\to H^2(M,\RR)\};
$$
\par {\rm (iii)} there exists a map $f: M\to K$ such that 
$$ 
[\omega]\in\Im\,\{f^*: H^2(K;\RR)\to H^2(M,\RR)\}.
$$
\end{prop}

\p  See \cite[Corollary 2.2]{RT}, cf. also \cite[Lemma 4.2]{LO}.
\hqed

\m As usual, we say that a closed form $\ga$ on $M$ is {\it integral} if the 
cohomology class $[\ga]\in H^*(M;\RR)$ comes from integral cohomology, i.e. 
that 
$\la \ga, a\ra \in \ZZ$ for all $a\in H_*(M)$.

\begin{prop}\label{integral}
Let $M$ be a symplectically aspherical manifold such that $H_2(M;\RR)$ is a 
finite dimensional vector space over $\RR$. Then $M$ admits an integral 
symplectically aspherical form. In particular, it holds if $M$ is a closed 
symplectically aspherical manifold.
\end{prop}

\p We represent $H_2(M;\RR)$ as a direct sum 
$$
H_2(M;\RR)=\Im(h(\pi_2(M))\otimes \RR)\oplus W
$$ 
where $h$ is the Hurewicz homomorphism. 

Let $\{a_1, \ldots , a_k\}$ be a family in $H_{2}(M)$ which yields an 
$\RR$-basis of $H_2(M;\RR)$. We can and shall assume that $W$ is spanned by 
$\{a_1, \ldots, a_l\}$, $l\le k$. Let $\gs_1, \ldots, \gs_k$ be closed 2-forms 
sush that $\la [\gs_i], a_j\ra=\delta_{ij}$. 

Choose a symplectically aspherical form $\omega$ on $M$ and consider the form 
$$
\eta=\omega + \gl_1\gs_1+\cdots + \gl_l\gs_l.
$$
It is clear that $\eta$ is a symplectic form for $\gl_i$ small enough and that 
$\eta|_{\pi_2(M)}=0$. Furthermore, we can always choose $\gl_i$ such that $\la 
\eta, a\ra\in \QQ$ for all $a\in H_2(M)$. Now, the form $N\eta$ is integral 
for 
$N\in \ZZ$ large enough.    
\hqed

\m Now we fix more properties of closed symplectically aspherical manifolds. 

\begin{cory}\label{fundclass}
Let $\omega$ be a symplectically aspherical form on a closed manifold 
$M^{2n}$, 
let $[M]\in H_{2n}(M;\RR)$ be the fundamental class,  and let $f: M \to 
K=K(\pi_1(M),1)$ be a map as in Proposition {\rm \ref{aspher}(iii)}. Then 
$f_*[M]\ne 0 \in H_{2n}(K;\RR)$. 
\end{cory}

\p Since $[\omega]=f^*(a)$ for some $a\in H^{2}(K;\RR)$, we conclude that
$$
0\ne \la \omega^n, [M]\ra= \la f^*(a^n), [M]\ra =\la a^n, f_*[M]\ra,
$$
and thus $f_*[M]\ne 0$. 
\hqed

\begin{cory}
The fundamental group of a closed symplectically aspherical manifold is 
infinite.
\end{cory}

\p Suppose that there exists a closed symplectically aspherical manifold $M$ 
with finite fundamental group. Passing to the (finite) universal covering, 
we can assume that $M$ is symply connected. Then the  Hurewicz 
homomorphism $\pi_2(M) \to H_2(M)$ is an isomorphism. On the other hand, 
$[\omega]\ne 0 \in H^2(M;\RR)$ since $M$ is closed, and so $\la \omega, 
a\ra\ne 
0$ for some $a\in H_2(M)=\pi_2(M)$. But this contradicts symplectic 
asphericity. 
\hqed

\part{Basic examples and constructions of symplectically aspherical manifolds}

\section{Homogeneous spaces}\label{sec-homogen}

Now we describe some classes of symplectically aspherical groups. For
this purpose, we recall several notions. 

A {\it lattice} in a Lie group $G$ is a discrete subgroup $\pi\subset G$. A
lattice $\pi$ in $G$ is called {\it uniform} if $G/\pi$ is compact. 

\begin{df} \rm A Lie group $G$ is called {\it completely solvable},
if any adjoint linear operator $ad\,V: \mathfrak g\to\mathfrak g$ of the Lie 
algebra $\mathfrak g$ of $G$ has only real eigenvalues.
\end{df}

It is well known that every completely solvable Lie group is solvable,
\cite{VGS} .

\begin{lemma}\label{compl-solv} 
Let $G$ be a completely solvable simply connected Lie group, and let $\pi 
\subset G$ be a uniform lattice in $G$. Suppose that the closed manifold 
$M=G/pi$ is c-symplectic. Then $M$ admits a symplectic structure. Finally, 
$M=K(\pi_1(M),1)$. In particular, $\pi_2(M)=0$, and therefore $M$ is 
symplectically aspherical. 
\end{lemma}

\p Since $M$ is c-symplectic, there exists a cohomology class $\alpha\in 
H^{2}(M,\RR)$ 
such that $\alpha^n\not=0\in H^{2n}(M;\RR)$. Since  $G$ is completely
solvable then, by the Hattori theorem \cite{Ha}, there is an isomorphism 
$$ 
H^*(M;\RR)\cong H^*(\Lambda\mathfrak g^*,\delta),
$$ 
 where $(\Lambda\mathfrak g^*, \delta)$ denotes the standard
Chevalley--Eilenberg complex for the Lie algebra $\mathfrak g$.  Therefore 
$\alpha$
can be represented by a closed differential 2-form $\omega$ whose pullback
$\widetilde \omega$ to $G$ is a left-invariant form.
Furthermore, $\widetilde\omega$ is non-degenerate since it is left-invariant 
and 
$[\widetilde\omega]^n\not=0$ on $H^{2n}(\Lambda\mathfrak g^*,\delta)$.
Hence, $\omega$ is non-degenerate, and so $(M, \omega)$ is a symplectic
manifold.  

Finally, $G$ is solvable, and hence $G$ is diffeomerphic to the Euclidean 
space. 
So, $\pi_k(G)=0$ for $k>1$, and thus $M=K(\pi_1(M),1)$ since $G$ is the 
universal cover of $M$. 
\hqed

\m Obviously, the class of completely solvable Lie groups contains all
nilpotent Lie groups. So, every c-symplectic nilmanifold $M$ 
is symplectically aspherical. (Recall that a nilmanifold is defined to be a 
manifold of the form $G/\pi$ where $\pi$ is a uniform lattice in a simply 
connected nilpotent Lie group $G$.)

\begin{cory}\label{4nil-aspher}
Every orientable $4$-dimensional nilmanifold $M$ is symplectically aspherical.
\end{cory}

\p It suffices to prove that $M$ is c-symplectic. It is well known that 
$H^2(M;\RR)\ne 0$. (You can use 
classification of 4-dimensional nilmanifolds, \cite{VGS}, or notice that 
$b_1(M)>1$ while 
$\chi(M)=0$.)   
Take any $a\in H^2(M;\RR), a\ne 0$. Then, by Poincar\'e duality, there exists 
$b\in 
H^2(M;\RR)$ with $ab\ne 0$. Since $ab=ba$, we have
$$
(a+b)^2=a^2+2ab+b^2,
$$
and so at least one of elements $a^2, b^2$ or $(a+b)^2$ must be non-zero. 
Thus, 
$M$ is a c-symplectic manifold.
\hqed

\begin{example}\label{ex-kahlerman}
 Here we show how to construct symplectically aspherical closed K\"ahler 
manifolds. Tori, products of complex curves and hyperplane sections of 
these manifolds give us examples of symplectically aspherical algebraic (and 
therefore K\"ahler) manifolds. Another source of such examples is based on the 
following construction. Let $G$ be a semisimple simply-connected Lie group of 
non-compact type, and let $K$ be a maximal compact connected subgroup of $G$. If 
the homogeneous space $G/K$ is a  symmetric Hermitian space then $G/K$ turns out to be a K\"ahler manifold with the invariant K\"ahler metric. All such pairs 
$(G,K)$ are listed in \cite[Ch. IX]{He}. Moreover, every such group $G$ contains 
a uniform lattice $\pi$, \cite{VGS}. Hence, $M:=\pi\backslash G/K$ is a closed 
K\"ahler manifold with the fundamental group $\pi$. Finally, $M=K(\pi_1(M),1$ 
since $G/K$ is diffeomorphic to Euclidean space, and thus $M$ is symplectically 
aspherical.
\end{example} 

\section{Gompf symplectic sum and symplectic asphericity}\label{sec-sum}

Here we discuss situations when the Gompf symplectic sum turns out to be a 
symplectically aspherical manifold.

\m First, recall the construction of the connected sum of two manifolds 
along a submanifold, with the aim to emphasize the symplectic version of this
construction, \cite{G1} . 

\m Let $M_1^n, M_2^n$ and $N^{n-2}$ be smooth closed oriented manifolds (not 
necessarily connected), of dimensions $n$ and $n-2$, respectively. Assume that 
we are given two embeddings $j_1: N\to M_1$ and $j_2: N\to M_2$, with the 
normal 
bundles $\nu_1$ and $\nu_2$, respectively, such that their Euler classes 
differ 
only by sign: 
$e(\nu_1)=-e(\nu_2)$. It turns out to be that there exists an 
orientation-reversing bundle isomorphism $\alpha: \nu_1 \to \nu_2$.  Let $V_i$ 
denote a tubular neighborhood of $j_i(N)$, which we identify with the total 
space of $\nu_i$. Then $\alpha$ yields a diffeomorphism $\psi: V_1\to V_2$,
which maps $j_1(N)$ to $j_2(N)$. Then $\psi$ determines an 
orientation-preserving diffeomorphism
$$
\varphi: (V_1-j_1(N))\to (V_2-j_2(N)), \,\,\,\varphi=\theta\circ\psi,
$$
where 
$$
\theta(p,v)=(p,\frac v {||v||^2})
$$ is a diffeomorphism which turns each
punctured normal fiber inside out.

\begin{df}\label{df-sum}
 \rm Let $M_1\cup_{\psi}M_2$ denote the smooth, closed oriented
manifold obtained from the disjoint union $M_1-(j_1(N))\sqcup M_2-(j_2(N))$ 
via 
gluing $V_1-j_1(N)$ and $V_2-j_2(N)$ by $\varphi$:
$$
M_1\cup_{\psi}M_2=M-(j_1(N)\cup j_2(N))/\simeq
$$
where $a\simeq b$ if and only if $b=\varphi(a),\,\,a\in V_1-j_1(N),b\in 
V_2-j_2(N).$
\end{df}

It was noted in \cite{G1} that there exists a cobordism $X$ between $M_1\sqcup 
M_2$ and 
$M_1\cup_{\psi}M_2$. It will be important for us to notice that the cobordism 
$X$ is 
obtained from
$(M_1\sqcup M_2)\times I$  ($I=[0,1]$) by identifying closed tubular 
neighborhoods of 
$j_1(N)\times 1$ and $j_2(N)\times 1$ by $\psi$ and rounding corners.

\m Now, we need the following observation. Every closed
$k$-form $\omega_M$ on $M$ for which $j_1^*\omega_M=j_2^*\omega_M$ induces a
cohomology class $[\Omega]\in H^k(X;\RR)$ and, hence, by restriction, a class
$[\omega]\in H^k(M_1\cup_{\psi}M_2;\RR)$. Note that $[\omega]=i^*[\Omega]$, 
where $i:M_1\cup_{\psi}M_2\to X$ is the canonical embedding.

In the sequel we will need the following result.

\begin{theorem}\label{sum}
Let $M_1, M_2, N$ and $j_i: N \to M_i, i=1,2$ be as in Definition 
$\ref{df-sum}$. Suppose in addition that $M_1, M_2$ and $N$ are symplectic 
manifolds and both embeddings $j_i: N \to M_i$ are symplectic. Then, for any 
choice of orientation reversing $\psi: V_1\cong V_2$, the manifold 
$M_1\cup_{\psi}M_2$ admits a canonical symplectic structure $\omega$, which is 
induced by the symplectic form on $M_1\sqcup M_2$ after a perturbation near
$j_2(N)$. More precisely, there is a unique isotopy class of symplectic forms
on $M_1\cup_{\psi}M_2$ $($independent of fiber isotopies of $\psi)$ that 
contains 
forms $\omega$ with the following characterization: the class $[\omega]\in
H^2(M_1\cup_{\psi}M_2;\RR)$ is the restriction of the class $[\Omega]\in 
H^2(X;\RR)$ canonically induced on the cobordism $X$ by the symplectic 
structure 
on 
$M_1\sqcup M_2$.
\end{theorem}

\p See Gompf \cite{G1}.
\hqed

\begin{theorem}\label{sum-sas}
Let $M_1,M_2,N$ and $j_i:N\to M_i$ be as in Theorem~$\ref{sum}$. Suppose in 
addition 
that $M_1, M_2$ and $N$ are connected.
If $M_1$ and $M_2$ are symplectically aspherical
and $(j_i)_*:\pi_1(N)\to \pi_1(M_i), i=1,2$ are monomorphisms,
then the symplectic sum 
$M_1\cup_{\psi }M_2$ is symplectically
aspherical.
\end{theorem}

\p Given a topological space $Y$, set $h(Y)=\Im(h:\pi_2(Y) \to H_2(Y)$ where 
$h$ is the Hurewicz homomorphism.
Let $X$ be the cobordism as above. Let $\omega$ be the form on $M_1\cup 
{\psi}M_2$ 
and $\Omega$ be the form on $X$. Since $\langle \omega, s\rangle= \langle 
\Omega, i_*s\rangle$, and $\Omega|_{M_i}=\omega_i$, it suffices to prove 
that 
$$
h(X)=\Im(h(M_1\times I)\to h(X))+\Im(h(M_2\times I)\to h(X)).
$$

\m
Because of general position arguments, every map $f:S^2\to X$ is homotopic 
to a smooth embedding. Moreover, we can always assume that $f(S^2)\cap W = 
f(S^2)\cap \Int W$ and that $f(S^2)$ intersects $W$ transvesally. So, 
consider an embedded sphere $S$ in $X$ which intesects $W$ transversely and 
does not meet $\pa W$. We must prove that the homology class of $S$ can be 
represented as as a disjoint union of spheroids in $X$, each of whish does not 
meet $W$.

\m Since $S\cap W$ must be a closed manifold, we conclude that $S\cap W$ is 
a disjoint union of circles, $S\cap W=C_1\sqcup \cdots \sqcup C_k$ where 
each $C_i$ is a circle. It is easy to see that $S\setminus (C_1\sqcup \cdots 
\sqcup C_k)$ contains at least one open 2-disc. Take one of such discs and 
denote it by $D$. Clearly, $D$ does not meet $W$. We denote the closure of 
$D$ in $S$ by $D_+$ and set $D_-=S\setminus D$. 

\m Clearly, $D_+\subset M_i\times I$ for a certain $i$, and we assume that 
$i=1$. Notice that the circle $C=\pa D_+$ is contractible in $M_1\times I$. 
So, it is contractible in $W$ because of monomorphicity of $(j_1)_*:\pi_1(N) 
\to \pi_1(M_1)$. Now we consider two different cases: $\dim W \ge 6$ and $\dim 
W=4$.

First, assume that $\dim W\ge 6$. Then there exists an 
embedded 2-disc $D_0$ in $W$ with $\pa D_0=C$. Moreover, we can even assume 
that $D_0$ does not meet $S\setminus C$.

\m Choose a trivialization $W\times [-1,1]$ of the normal bundle of $W$ in 
$X$. We can assume that $S\cap (D\times [-\eps, \eps])C\times [-\eps, \eps]$ 
for 
$\eps$ small enough. 
Now, set 
$$
S_1=(D_+\setminus (C\times [0,\eps)))\cup D_0\times \{\eps\}
$$
and 
$$
S_2=(D_-\setminus (C\times (-\eps,0]))\cup D_0\times \{-\eps\}.
$$ 
Clearly, $S_1$ and $S_2$ are embedded 2-spheres, and the homology class of 
$S$ is the sum of the homology classes of $S_1$ and $S_2$. 
Furthermore, $S_1$ sits in $M_1\times I$, while the intersection $S_2\cap W$ 
has less connected components then $S\cap W$ does. Now we leave $S_1$ as it is 
and perform the same procedure as above with $S_2$, and so on. Finally, we get 
the disjoint union of spheroids each of which does not meet $W$, and this 
disjoint union represents the same homology class as $S$.

\m Finally, if $\dim M =4$ then we can find not embedded but immersed disk 
$D_0$ 
in $W$, and we can assume the the interior of $D_0$ does not meet $S$ (since 
$\dim X=5$). Now we can do the same procedure as above, but $S_1$ and $S_2$ 
will 
be not embedded but immersed sphere. We can leave $S_1$ as it is, while we can 
perform a small perturbation of $S_2$ and get an embedded sphere $S'_2$ such 
that $S_2\cap W=S'_2\cap W$. Now we can finish the proof as in the case $\dim 
W\ge 6$.      
\hqed

Let $(X;A,B)$ be a $CW$-triad. We set $C=A\cap B$ and denote by $j_1: A \to 
X$, 
$j_2: B \to X$, $i_1: C \to A$, $i_2: C \to B$ the obvious inclusions.

\begin{prop}\label{mv}
 Fix any $k$ and any coefficient group. If the homomorphism $i_1^*: H^k(A) \to 
H^k(C)$ is injective and the homomorphism $i_1^*: H^{k-1}(A) \to H^{k-1}(C)$ 
is 
surjective, then the $j_2^*: H^k(X) \to H^k(B)$ is injective and the 
homomorphism $j_2^*: H^{k-1}(X) \to H^{k-1}(B)$ is surjective.
\end{prop}

\p The exactness of the sequence
$$
\CD
H^{k-1}(A) @>i_1^*>> H^{k-1}(C) \to H^{k}(A,C) @>i_1^*>> H^k(A) \to H^k(C)
\endCD
$$
implies that $H^{k}(A,C)=0$. So, because of the excision property, 
$H^{k}(X,B)= H^{k}(A,C)=0$. Now, the exactness of the sequence
$$
\CD
H^{k-1}(X) @>j_2^*>> H^{k-1}(B) \to H^{k}(X,B) @>i_2^*>> H^k(X) \to H^k(B)
\endCD
$$
implies the required claim on $j_2^*$.
\hqed
 
\m We say that a 2-dimensional cohomology class $a$ is {\it decomposable} if 
it can be represented as $a=\sum_i a_i a'_i$ where $a_i$ and $a'_i$ are 
1-dimensional classes. Notice that a symplectic form $\omega$ on a symplectic 
manifold is aspherical if its cohomology class $[\omega]$ is decomposable.

\begin{theorem} Let $M_1,M_2,N$ and $j_i: N \to M_i, i=1,2$ be as in Theorem 
$\ref{sum}$, and suppose that $j_1$ induces a surjection 
on the first cohomology group and an injection on the second cohomology group.  
Assume that $H^2(M_2;\RR)$ consists of decomposable elements.
Then the symplectic manifold $(M_1\cup_{\psi}M_2)$ is symplectically 
aspherical. 
\end{theorem}

\p Let $X$ be the cobordism described after Definition \ref{df-sum}. It 
suffices to prove that the cohomology class $[\Omega]$ of the form $\Omega$ on 
$X$ is decomposable. Notice that $X$ is homotopy equivalent to the space 
$Y=M_1\cup _N  M_2$. So, we have the triad $(Y; M_1, M_2)$ with $M_1 \cap 
M_2=N$. 
Since smooth  manifolds are triangulable, we can regard the above triad as a 
$CW$-triad. Let $j: M_2 \to Y$ be the inclusion. Because of the conditions of 
the Theorem, 
$$
j^*[\Omega]=\sum b_i b'_i, \quad b_i, b'_i \in H^1(M_2).
$$

By the Proposition \ref{mv}, the map $j^*:H^1(Y) \to H^1(M_2)$ is an 
epimorphism. So, 
there are $a_i, a_i\in H^1(Y)$ with $j^*(a_i)=b_i, j^*(a'_i)=b'_i$. So, 
$$
j^*([\Omega] - \sum a_i a'_i)=0.
$$ 
But, again by the Proposition \ref{mv}, the 
map $j^*:H^2(Y) \to H^2(M_2)$ is a monomorphism, and thus $[\Omega] - \sum a_i 
a'_i=0$.
\hqed

\m Another situation when the Gompf symplectic sum 
is symplectically aspherical is described in \cite[Proposition 3.1]{K}.

\section{Symplectic bundles and symplectic asphericity}\label{sec-bun}

We define a {\it symplectic bundle} to be a locally trivial bundle whose fiber 
is a closed symplectic manifold $F=(F,\omega)$ and with the structure group 
$\sympl(F,\omega)$ of symplectomorpisms of $F$. 
The goal of the section is to prove the following theorem.

\begin{theorem}\label{main-bundle}
Let
$$
\CD 
\xi=\{F @>{i}>> M @>{p}>> B\}
\endCD
$$
be a symplectic bundle. Suppose also that there exists a 
cohomology class $a\in H^2(M;\RR)$ such that $i^*a=[\omega_F]$. Then the 
following 
holds:

{\rm (i)} the group $\pi_1(M)$ is an
extension of $\pi_1(B)$ by $\pi_1(F)$, i.e. there is an exact sequence
$$
1\longrightarrow \pi_1(F)\longrightarrow \pi_1(M)\longrightarrow
\pi_1(B)\longrightarrow 1;
$$

\par {\rm (ii)}  if the form $\omega_F$ is symplectically aspherical and $B$ 
is 
closed symplectically aspherical manifold such that the Hurewicz homomorphism 
$h:\pi_2(B)\otimes \RR \to H_2(B;\RR)$ is monic, then $M$ is symplectically 
aspherical.
\end{theorem}

\m We preface the proof with some supporting claims.

\begin{lemma}\label{aspher-bundle}
Let
$$
\CD 
F @>{i}>> E @>{p}>> B
\endCD
$$
be a locally trivial bundle. Suppose that there is a class $a\in H^2(F;\RR)$ 
such that $i^*(a)|_{\pi_2(F)}=0$. Assume also that the Hurewicz homomorphism 
$h:\pi_2(B)\otimes \RR \to H_2(B;\RR)$ is monic. 
Then there exists a class $b\in H^2(E;\RR)$ with $i^*b=i^*a$ and 
$b|_{\pi_2(E)}=0$.
\end{lemma}

\p Consider the long exact sequence
$$
\CD
\cdots @>>> \pi_2(F) @>{i_*}>> \pi_2(E) @>{p_*}>>\pi_2(B) @>>> \cdots
\endCD
$$
and decompose $\pi_2(E)\otimes \RR=i_*(\pi_2(F)\otimes \RR)\oplus W$.

Choose a base $\{a_1,\ldots,a_k\}$ of $h(W)$ where $h: \pi_2(E)\otimes \RR \to 
H_2(E;\RR)$ 
is the Hurewicz homomorphism. Define 
$$
\mu_i=<a,a_i>,\quad i=1,\ldots ,k.
$$

Since the homomorphism $h:\pi_2(B)\otimes \RR \to H_2(B;\RR)$ is monic, the 
family $\{p_*(a_j)\}$ is linearly independent in $H_2(B;\RR)$. Let $c_i\in 
H^2(B;\RR)$ be elements dual to $p_*(a_i)\in H_2(B;\RR)$,
$$
<c_i,p_*(a_j)>=\delta_{ij}.
$$
We set
$$
b=a-\sum_{i=1}^k \mu_i p^*c_i.
$$

Clearly, $i^*b=i^*a$. So, $b|_{i_*(\pi_2(F))}=0$, and it remains to prove 
that $b|_{W}=0$, i.e, that $<b,a_i>=0$ for all $i$. But
\begin{eqnarray*}
<b,a_i> &=&\left\langle a-\sum_{j=1}^k\mu_j p^*c_j,a_i\right\rangle 
=<a,a_i>-\sum_{j=1}^k\mu_j <p^*c_j,a_i>\\
&=&\mu_i-\mu_i=0,
\end{eqnarray*}
and we are done.
\hqed

\m Let $G$ be a topological group. Given a map $\gf: S^{n-1}\to G$ and a 
$G$-space $F$, the clutching construction gives us a locally trivial 
bundle $\eta=\{F \to E \to S^n\}$ with the structure group $G$, \cite{St}. 
Choose a point $x_0\in F$ and define the evaluation map
$$
\ev : G \to F, \quad \ev(g)=g(x_0).
$$
If we regard $G$ as the pointed space with the unit as the base point, then 
$\ev$ turns into a pointed map. In particular, it induces a map $\ev_*: 
\pi_{n-1}(G) \to \pi_{n-1}(F)$. Now, in the exact sequence
$$
\CD
\pi_n(F) \to \pi_n(E) @>p>> \pi_n(S^n) @>\pa>>\pi_{n-1}(F) \to
\endCD
$$
we have 
\begin{equation}\label{eval}
\pa\ga=\ev_*([\gf])
\end{equation}
where $\ga \in \pi_n(S^n)=\ZZ$ is the generator and $[\gf]\in \pi_{n-1}(G)$ is 
the homotopy class of $\gf$ , cf. \cite[\S 18]{St}.

\begin{prop} [{\cite[Section 1.2]{LMP}}]
\label{lmp}
Let $M=(M,\omega)$ be a closed symplectic manifold, and let $\ham(M)$ be the 
group of hamiltonian symplectomorphisms of $M$. Then $\ev_*: \pi_1(\ham(M)) 
\to 
\pi_1(M)$ is the zero homomorphism.
\hqed
\end{prop}

\begin{cory}\label{zero-map}
Let $F\to E\to S^2$ be a symplectic bundle whose structure group reduces to 
$\ham(F)$. Then the boundary map $\pa :\pi_2(S^2)\to \pi_1(F)$ in the homotopy 
exact sequence of the bundle is the zero map.  
\end{cory}

\p It follows from Propsition~\ref{lmp} and equality \eqref{eval}
\hqed

\begin{theorem}[{\cite[Theorem 6.36]{MS}}]
\label{red-ham}
Let 
$$
\CD
F @>i>> E @>p>> B.
\endCD
$$  
be a symplectic bundle over a simply connected base $B$. Assume that there 
exists a class $a\in H^2(E)$ with $i^*(a)=[\omega_F]$. Then the structure 
group of the bundle can be reduced to $\ham(F)$.
\hqed
\end{theorem}

\m Now we are ready to prove Theorem \ref{main-bundle}.

\p
(i) It suffices to prove that, in the 
exact sequence 
$$
\CD
\cdots \to \pi_2(M) \to \pi_2(B) @>\pa >> \pi_1(F) \to \cdots
\endCD
$$
of the bundle $\xi$, we have $\pa \gb=0$ for every $\gb\in \pi_2(B)$. We 
represent $\gb$ by a map $f: S^2 \to B$ and consider the bundle 
$$
f^*\xi=\{F \to E \to S^2\}
$$
induced from the bundle $\xi$ by $f$. Then we have 
a commutative diagram 
$$
\CD
 F  @>i'>> E  @> >> S^2 \\
 @| @VgVV @VVfV \\
 F  @>i>> M  @>>> B 
\endCD
$$ 
of bundles. Clearly, $(i')^*g^*(a)=[\omega_F]$, and so, by \ref{red-ham}, the 
structure group of the top bundle reduces to $\ham(F)$. Now, we have a 
commutative diagram
$$
\CD
\cdots @>>> \pi_2(E) @>>> \pi_2(S^2) @>\pa' >> \pi_1(F) @>>> \cdots\\
@. @VVV @Vf_*VV @| @.\\
\cdots @>>> \pi_2(M) @>>> \pi_2(B) @>\pa >> \pi_1(F) @>>> \cdots
\endCD
$$ 
of exact sequences. By \ref{zero-map}, $\pa'\ga =0$. Now, 
$\pa \gb=\pa f_*\ga =\pa'\ga =0$, as asserted.

\par (ii)  Let $\omega_F$ and $\omega_B$ be the symplectically aspherical 
forms 
on $F$ and $B$, respectively. We apply Lemma 
\ref{aspher-bundle} to the bundle $\xi$ and see that there exists a class 
$b\in 
H^2(M;\RR)$ such that $i^*b=[\omega_F]$ and $b|_{\pi_2(M)} = 0$. Because of 
a Thurston Theorem, see e.g. \cite[Theorem 6.3]{MS}, $M$ possesses a 
symplectic 
form 
$\omega$ with $[\omega]=b+Kp^*[\omega_B]$ with $K\in \RR$. Clearly, this form 
$\omega$ is symplectically aspherical.  
\hqed

\part{Fundamental groups of closed symplectically aspherical manifolds}

Let us call a group symplectically aspherical if it can be realized as a 
fundamental group of a closed symplectically aspherical manifold. Because of 
Proposition~\ref{induced}, any subgroup of finite index of a symplectically 
aspherical group is symplectically aspherical.

\m From here and to the end of the paper, ``symplectic manifold''  means 
``closed symplectic manifolds'' unless some other is said explicitly.

\section{Dimension phenomena}\label{phen}

The following theorem is a symplectic analog of the Lefschetz Theorem on 
Hyperplane Sections.

\begin{theorem}
\label{donaldson}
 Let $(M^{2n},\omega)$ be a symplectic manifold with an integral 
symplectic form $\omega$, and let
$h\in H^2(M)$ be an integral lift of $[\omega]$. Then for $N$ large enough the
Poincar\'e dual of $Nh$, in $H_{2n-2}(M)$, can be realized by a symplectic
submanifold $V^{2n-2}$ of $M^{2n}$. Moreover, we can choose $V$ such that the
inclusion $i: V \hookrightarrow M$ is an $(n-2)$-equivalence, i.e. the
homomorphism $i_*: \pi_k(V) \to \pi_k(M)$ is an isomorphism for $k\le n-2$ and
an epimorphism for $k=n-1$. 
\end{theorem}

\p See \cite[Theorem 1 and Proposition 39]{D}.
\hqed

\begin{prop}
\label{down} 
Suppose that a group $\pi$ can be realized as the
fundamental group of a symplectically aspherical manifold $(M^{2n},\omega)$
with $n\ge 3$. Then $\pi$ can be realized as the fundamental group of a
$(2n-2)$-dimensional symplectically aspherical manifold. 
\end{prop}

\p  Because of Proposition~\ref{integral}, we equip $M$ with an integral 
symplectically aspherical form.  According to the Theorem 
\ref{donaldson}, there exists a symplectic submanifold $V^{2n-2}$ of $M$ such 
that the inclusion $i: V \hookrightarrow M$ is an $(n-2)$-equivalence. In 
particular, $\pi_1(V)=\pi_1(M)$.  Clearly, $V$ is symplectically aspherical 
since $M$ is, and the result follows.
\hqed

\begin{cory}
\label{four}
Suppose that a group $\pi$ can be realized as the
fundamental group of a symplectically aspherical manifold $(M^{2n},\omega)$
with $n\ge 3$. Then $\pi$ can be realized as the fundamental group of a
$4$-dimensional symplectically aspherical manifold. 
\hqed
\end{cory}

\m So, we can decrease the dimension of the symplectically aspherical manifold 
with a prescribed fundamental group. However, we are not always able to 
increase the dimension, as the following proposition shows.

\begin{prop} \label{dim}
Let $\pi$ be a group such that $H^{2k}(\pi;\RR)=0$. Suppose that $\pi$ 
can be realized  as the fundamental group of a symplectically aspherical 
manifold $M^{2n}$. Then $n\le k$. 
\end{prop}

\p Choose a symplectically aspherical form $\omega$ on $M$. Because of 
\ref{aspher}, there exists a map $f: M \to K(\pi,1)$ such that 
$$ 
[\omega]\in \Im \{f^*: H^2(\pi;\RR)=H^*(K(\pi,1);\RR) \to H^2(M;\RR)\}.
$$ 
Since $[\omega]^{2n}\ne 0$, we conclude that $n \le k$.
\hqed 

\begin{cory} The group $\ZZ^{m}$ cannot be realized as the
fundamental group of a symplectically aspherical manifold of dimension $2k$
with $2k>m$.
\hqed
\end{cory}

\m Notice that $\ZZ^{2n}$ is the fundamental group of the torus $T^{2n}$. 
Since
$\pi_2(T^{2n})=0$,  $\ZZ^{2n}$ can be realized as the fundamental group of a
symplectically aspherical manifold of dimension $2k$ with $2\le k \le n$.

\begin{remark} \rm Because of Propositions \ref{four} and \ref{dim}, it makes 
sense to introduce the following invariant of symplectically aspherical 
groups. 
Namely, given a symplectically aspherical group $\pi$, we define $\nu(\pi)$ to 
be the largest $n$ such that $\pi$ can be realized as the fundamental group of 
a closed symplectically aspherical manifold $M^{2n}$. For example, 
$\nu(\ZZ^{2n})=n$. Furthermore, if $\pi$ is the fundamental group of a closed 
orientable surface then $\nu(\pi)=1$, and if $G$ is the direct product of $n$ 
such groups then $\nu(G)=n$.
\end{remark}

\section{Two classes of symplectically aspherical groups}\label{classes} 

Let $\AAA$ be the class of symplectically aspherical groups which can be
realized as the fundamental groups of symplectically aspherical manifolds with
trivial $\pi_2$, and let $\BBB$ be the class of symplectically aspherical
groups which can be realized as the fundamental groups of symplectically
aspherical manifolds with non-trivial $\pi_2$. In this section we want to
investigate the relation between the classes $\AAA$ and $\BBB$. First, some
trivial remarks.

\begin{enumerate}
\item If $\pi \in \BBB$ and $\tau$ is simplectically aspherical, then $\pi
\times \tau \in \BBB$.
\item Let $G$ be the fundamental  group of a closed orientable surface. Then 
$G
\notin \BBB$ (by Proposition \ref{dim}). 
\end{enumerate}

\m
According to the results of Section \ref{phen}, it seems 
reasonable to introduce the classes $\AAA_{2n}$ an 
$\BBB_{2n}$ as follows. The group $\pi$ belongs to $\AAA_{2n}$ if $\pi$ can be 
realized as the fundamental group of a symplectically
aspherical $2n$-dimensional manifold $M$ with $\pi_2(M)=0$. Similarly, the group 
$\pi$ belongs to $\BBB_{2n}$ if $\pi$ can be realized as the fundamental group 
of a symplectically
aspherical $2n$-dimensional manifold $M$ with $\pi_2(M)\ne 0$. Because of 
what is done above, we have the following Proposition.

\begin{prop}\label{chain}
$$
\AAA_{2n+2} \subset \AAA_{2n} \subset \cdots \subset \AAA_6\subset \AAA_4\cup 
\BBB_4, \quad  \BBB_{2n+2}
\subset \BBB_{2n} \subset \cdots \subset \BBB_4.
$$
\end{prop}

\p Consider the $(n-2)$-equivalence $i: V^{2n-2} \to M^{2n}$ from the proof of 
Theorem \ref{down}. The homomorphism $i_*: \pi_2(V) \to \pi_2(M)$ is an 
isomorphism for $n\ge 4$ and an epimorphism for $n=3$. Therefore the result 
holds.
\hqed

Recall the following theorem of Hopf.

\begin{theorem}
\label{hopf}
 Let $X$ be a connected $CW$-space with
$\pi_1(X)=\pi$ and $\pi_i(X)=0$ for $i<n$. Then there is an exact sequence
$$
\CD
\pi_n(X) @>h>> H_n(X) @>>> H_n(\pi) @>>> 0.
\endCD
$$   
\end{theorem}

\p See \cite[Theorem II.5.2]{B}.
\hqed

We denote by $b_i(X)$ the $i$-th Betti number of (the group or the space) $X$.

\begin{theorem}
\label{sasgroup}
Let $\pi $ and $\tau $ be symplectically aspherical groups.
Then the following statements are true:

\begin{enumerate}
     \item[(i)] If $\pi \in {\AAA}_4$, then $b_1(\pi )\geq b_3(\pi)$.

     \item[(ii)] If $b_1(\pi ) < b_3(\pi )$, then $\pi \in {\BBB}$.

     \item[(iii)] Let $\pi $ be fundamental group of a closed symplectically
           aspherical manifold $(M,\omega )$ with $\pi _3(M)=0$.
           If $b_1(M) < b_3(M)$ then $\pi \in {\BBB}$.

     \item[(iv)] If $ \max\{b_3(\pi ),b_3(\tau )\} \geq 1$,
           then $\pi \times \tau \in {\BBB}$.

     \item[(v)] If $\max\{b_2(\pi ),b_2(\tau )\} \geq 2$ and
           $\max\{b_1(\pi ),b_1(\tau )\} \geq 1$ 
           then $\pi \times \tau \in {\BBB}$.

     \item[(vi)] $\pi \times \mathbb Z^4 \in {\BBB}$.
   \end{enumerate}
\end{theorem}

\p (i) Because of the Hopf Theorem \ref{hopf}, there is an epimorphism $H_3(M) 
\to
H_3(\pi)$. So, $b_1(\pi)=b_1(M)=b_3(M)\ge b_3(\pi)$.

\par (ii)
Because of Corollary \ref{four}, $\pi$ can be realized as the fundamental 
group of a
$4$-dimensional symplectically aspherical manifold. Now the result follows 
from (i).

\par (iii) If $\pi_2(M)\not=0$ then we are done. So, suppose that 
$\pi_2(M)=0$. Since
$\pi_3(M)=0$, the Hopf exact sequence from Theorem \ref{hopf} for $n=3$ yields 
an isomorphism
$H_3(M) \cong H_3(\pi)$. So, 
$$ b_1(\pi)=b_1(M) < b_3(M) = b_3(\pi)
$$ and the result follows from (ii). 

\par (iv), (v) Notice that $b_2(G)>0$ for every symplectically aspherical 
group $G$. Now,
$b_1(\pi\times \tau)=b_1(\pi) + b_1(\tau)$, while (by the K\"unneth formula)
$$ b_3(\pi \times \tau) = b_3(\tau) + b_1(\pi)b_2(\tau) + b_2(\pi)b_1(\tau)+
b_3(\pi).
$$ Now, each of the conditions in (iv), (v) implies that 
$$ b_3(\pi \times \tau) > b_1(\pi \times \tau),
$$ and  the result follows from (iii).

\par (vi) This follows from (v) directly.
\hqed

\begin{remque}\rm 
\m 1. Notice that $\BBB_2=\emptyset \ne \AAA_2$.  

\m 2. We have $\AAA_6 \cap \BBB_6 \ne \emptyset$ since $\ZZ^{8} \in \AAA_6\cap 
\BBB_6$. Indeed,  $\ZZ^{8}\in
\AAA_{8} \subset \AAA_6$. On the other hand, $\ZZ^6 =\pi_1(T^6)$ and therefore 
$\ZZ^{6}\in \BBB_{4}$ by Theorem \ref{sasgroup}(i). Thus, $\ZZ^8 =\ZZ^6 \times 
\ZZ^2 
\in \BBB_6$. 

\m 3. Similalry, $\AAA_{2n} \cap \BBB_{2n} \ne \emptyset$ for
$n\ge 3$. 

\m 4. We don't know whether $\AAA_4 \cap \BBB_4 \ne \emptyset$. In
particular, is it true that $\ZZ^4 \in \BBB_4$?    

\m 5. Is $\ZZ^{2n+1}$ symplectically aspherical if $n>1$? ($\ZZ$ and $\ZZ^3$ 
are
not by Proposition \ref{dim}.) If the answer is negative, the proof should be 
delicate because the answer is positive at c-symplectical level, see 
Proposition \ref{zz} below.

\m 6. Generally, is it true that $\BBB \subset \AAA$?
\end{remque}

\begin{prop}\label{zz}
For every $n>3$ there exists a closed manifold $N^{2n}$ and a cohomology class 
$a\in 
H^2(N;\RR)$ such that $a^n\ne 0$, $\pi_1(N)=\ZZ^{2n+1}$ and $a|_{\pi_2(N)}=0$.
\end{prop}

\p Take the torus $T^{2n+2}$ and consider its hyperplane section as in Theorem
\ref{donaldson}. Then we get a 2n-dimensional symplectically aspherical 
manifold 
$(M,\omega)$ with the fundamental group $\ZZ^{2n+2}$. Then, by Proposition 
\ref{aspher}, there exists a  map $f: M\to T^{2n+2}$ which induces an 
isomorphism of fundamental groups and such that $f^*(u^{n})=[\omega]^n\ne 0$ 
for 
some $u\in H^2(T^{2n+2};\RR)$. So, there are cohomology classes $x_i\in 
H^1(T^{2n+2};\RR), i=1,2, \ldots 2n$ such that $f^*(x_1 \cdots x_{2n})\ne 0$. 
This implies, in  turn, that there exists a map $g: M \to T^{2n}$ with 
$g^*[\omega_T]^n\ne 0$. Here $\omega_T$ is the symplectic form on $T^{2n}$. In 
particular, the degree of $g$ is 
non-zero. 

\m
Consider the induced homomorphism 
$$
g_*:\ZZ^{2n+2}=\pi_1(M) \to  \pi_1(T^{2n})=\ZZ^{2n}
$$ 
and take any $b\in \Ker g_*$. Let $A$ be the subgroup generated by $b$. Then 
$$
\ZZ^{2n+2}/A \cong \ZZ^{2n+1} \oplus F
$$ 
where $F$ is a finite abelian group.

\m
Now, we represent $b$ by an embedded circle $S$ and perform the surgery of $g$ 
along $S$. Then we get a map $h: P \to T^{2n}$ which is bordant to $g$, and 
therefore $h$ has non-zero degree.  So, $h^*[\omega]^n\ne 0$, and thus $P$ is 
c-symplectic. Furthermore, $\pi_1(P)=\ZZ^{2n+1} \oplus F$. Now, consider a 
finite cover $\gf: N \to P$ which corresponds the direct summand $Z^{2n+1}$ of 
$\pi_1(P)$ and put $a=\gf^*(h^*\omega)$. Then $a^n\ne 0$ since 
$h^*([\omega]^n)\ne 0$. Finally, the class $a$ vanishes on the image of the 
Hurewicz map since $[\omega]$ does.
\hqed

\section{Some results about realization} 
 
Now we describe some classes of symplectically aspherical groups. Let us call 
a group c-symplectic if $K(\pi, 1)$ can be realized as a closed c-symplectic 
manifold. 

\begin{prop}\label{c-sympl} 
If $\pi$ is a uniform lattice in a simply-connected completely solvable Lie 
group $G$ of dimension $2n$ and $\pi$ is 
c-symplectic, then $\pi \in
\AAA_{2n}$. In particular, $\pi$ is symplectically aspherical. 
\end{prop}

\p Since every $G$ is difeomorphic to $\RR^{2n}$, This is a direct corollary 
of Lemma~\ref{compl-solv}, because $G/\pi$ is the Eilenberg--Mac Lane space 
$K(\pi,1)$.
\hqed

\m Let $\pi$ be a polycylic group. Let $\alpha\in\Aut(\pi)$. There exists a
subnormal series $\pi=\pi_n\supset \pi_{n-1}\supset ... \supset \pi_0$ such
that $\alpha(\pi_i)\subset\pi_i$, \cite{Gb}. (Here subnormality means that
$\pi_i$ is normal in $\pi_{i+1}$ and $F_i=\pi_{i+1}/\pi_i$ are finitely
generated abelian groups.) Hence $\alpha$ induces automorphisms $\alpha_i\in
\Aut(F_i \otimes\CC)=GL(k_i,\CC)$. One can easily check that the set of
eigenvalues of all operators $\alpha_i$ does not depend on the choice of a
subnormal series. We call the elements of this set eigenvalues of $\alpha$. 

\begin{df} A polycyclic group $\pi$ is called a group of type
($R$), if for all $\gamma\in\pi$ all eigenvalues of the inner automorphism 
$\rm{Int}(\gamma)$ are real and positive. 
\end{df}

\begin{theorem} A group  $\pi$ is isomorphic to a uniform lattice in a
completely solvable simply-connected Lie group if and only if $\pi$ is of type
$(R)$.
\end{theorem}

\p See Gorbatsevich \cite{Gb}.
\hqed

\begin{cory} Every c-symplectic group $\pi$ of type $(R)$ belongs to $\AAA$. 
Furthermore, $\pi \in \BBB$ if
$b_1(\pi)<b_3(\pi)$.
\hqed
\end{cory}

\m It is well known that every finitely
generated torsion free nilpotent group is of type ($R$), \cite{M, VGS}, and so 
every such group belongs to $\AAA$ provided it is c-symplectic. Now we show 
that 
some of these groups really belong to $\BBB$.

\begin{cory} 
\label{six}
The fundamental group of any $6$-dimensional symplectic
nilmanifold is a symplectically aspherical group of class $\BBB$.
\end{cory}

\p All 6-dimensional symplectic nilmanifolds are classified (see \cite{Sa,
IRTU}). In particular, the first and second Betti numbers of each of 34 such
manifolds can be found in the corresponding tables in these papers. Note that
since the Euler characteristic of any nilmanifold is zero, we get the 
following
relation for the Betti numbers: $2-2b_1+2b_2-b_3=0$. Hence $b_1<b_3$ is the
same as 
$$
2+2b_2>3b_1.
$$
One can check that each of the symplectic nilmanifolds from the tables
\cite{Sa, IRTU} satisfies this inequality.
\hqed

\m 
 Notice that one can also get groups of type $(R)$ which are solvable but
non-nilpotent.
 For example, consider the following simply-connected completely solvable Lie
group $G$ consisting of matrices
$$
\begin{pmatrix} e^t & 0 & xe^t & 0 & 0 & y_1\\ 0 & e^{-t} & 0 & xe^{-t} & 0 & 
y_2\\ 0
& 0 & e^t & 0 & 0 & z_1\\ 0 & 0 & 0 & e^{-t} & 0 & z_2\\ 0 & 0 &  0 & 0 & 1& 
t\\
0 & 0 & 0 &  0 & 0 & 1
\end{pmatrix} .
$$ It is shown in \cite{FLS}, that this group contains a uniform lattice
$\pi$, and that the compact solvmanifold $M:=G/\pi$ has 
$b_1(M)=2<b_3(M)=4$.
Thus, $\pi \in \BBB$.

\begin{example} In Example \ref{ex-kahlerman} we constructed a (closed) 
symplectically aspherical K\"ahler manifolds of the form $\pi\backslash G/K$ 
with the fundamental group $\pi$ and trivial $\pi_2$. So, $\pi\in \AAA$. Now we 
give an example of $\pi$ as above  with $\pi \in \BBB$.

Let $\chi(\pi)$ denote the Euler characteristic of $\pi$. It was shown in
\cite[Theorem 7.9]{VGS} that $\chi(\pi)\not=0$ if and only if
$\text{rank}(G)=\text{rank}(K)$, and in the latter case one has also the sign
of $\chi(\pi)$  equal to $(-1)^n$, where $n= {1/2}\dim G/K$. Now,
consider $G=Sp(2,\RR)$ and $K=U(2)$.  Then $G/K$ is a 6-dimensional Hermitian
symmetric space of non-compact type, and therefore $\chi(\pi)<0$.
Furthermore, $b_1(\pi)=0$, see \cite[Theorem 7.1]{VGS}. On the other hand,
$\chi(\pi)=2-2b_1+ 2b_2-b_3=2+2b_2-b_3<0$, which implies $b_3(\pi)>0$.
Thus, $\pi \in \BBB$.
\end{example} 

\begin{example}
The group theoretic constructions which might occur as
symplectically aspherical groups.
\begin{enumerate}
\item
Obviously, the product of two symplectically aspherical
groups is symplectically aspherical.
\item
Hopefully, other certain extensions
of symplectically aspherical groups are again symplectically
aspherical, in view of Theorem~\ref{main-bundle}.
\item
Amalgamated products can occur due to Theorem~\ref{sum-sas}
In particular, if $F_1$, $F_2$ are closed orientetd
surfaces and $\pi $ is symplectically aspherical group,
then the group

$$
\left (\pi_1(F_1)\times \pi \right ) *_{\pi }
\left (\pi \times \pi_1(F_2)\right )
$$
\noindent
is symplectically aspherical by the Seifert--van Kampen Theorem .
\end{enumerate}
\end{example}

\section{Nilpotent groups in $\AAA_4$} 

In this section we describe all nilpotent groups which can be realized as the 
fundamental groups of symplectic manifolds $M$ with $\pi_2(M)=0$. Here we use 
some ideas from \cite{R1}.

\m Let $\pi$ be a finitely presented group, and let $X$ be a $CW$-space with 
$\pi_1(X)=\pi$ and finite 2-skeleton. Let $\wt X$ be the universal covering 
space of $X$, and let $H^1_c(\wt X)$ be the 1-dimensional cohomology with 
compact supports of $\wt X$. 

\begin{propdef}\label{propdef}
The group $H^1_c(\wt X)$ depends on the group $\pi$ only. We denote it by 
$H^1_c(\pi)$ and call the 1-dimensional cohomology with compact supports of 
$\pi$. 
\end{propdef}

\p Consider two spaces $X_1$ and $X_2$ as the above described space $X$. 
First, 
assume that both $X_1$ and $X_2$ are $K(\pi,1)$'s. Consider homotopy 
equivalences $f: X_1 \to X_2$ and $g: X_2 \to X_1$ with $gf \simeq 1_{X_1}$ 
and 
$fg \simeq 1_{X_2}$. We can assume that $f(X_1^{(2)})\subset X_2^{(2)}$ and 
$g(X_2^{(2)})\subset X_1^{(2)}$. Moreover, the homotopies $H: gf\simeq 1$ and 
$H': fg \simeq 1$ can be chosen so that $H(X_1^{(1)} \times I) \subset 
X_1^{(2)}$ 
and $H'(X_2^{(1)} \times I) \subset X_2^{(2)}$. 

\m Passing to the universal coverings, we get the homotopy equivalences $\wt 
f: 
\wt X_1 \to \wt X_2$ and homotopies $\wt H: \wt X_1 \times I \to \wt X_1$ and 
$\wt H': \wt X_2 \times I \to \wt X_2$. Clearly, the maps 
$$
\wt f|_{X_1^{(2)}}: X_1^{(2)} \to X_2,\quad \wt g|_{X_2^{(2)} }:X_2^{(2)} \to 
X_1
$$
and the homotopies
$$
\wt H|_{X_1^{(1)} \times I}: {X_1^{(1)} \times I} \to X_1, \quad \wt 
H'|_{X_2^{(1)} \times I}: {X_2^{(1)} \times I} \to X_2
$$
are proper maps. Therefore $\wt f$ induces an isomorphism $H^1_c(\wt X_2) \to 
H^1_c(\wt X_1)$.

\m Now we consider an arbitrary space $X$ as above. We attach to $X$ cells of 
dimension $\ge 3$ and get an embedding $X \subset Y$ where $Y=K(\pi,1)$. Since 
$X^{(2)}=Y^{(2)}$, we conclude that $H^1_c(\wt X) = H^1_c(\wt Y)$. This 
completes the proof.
\hqed

\begin{remark}\rm
Certainly, the group  $H^1_c(\pi)$ admits a purely algebraic description in 
terms of the group $\pi$, cf \cite{R1, N}. However, the description from 
\ref{propdef} is enough for our goals.
\end{remark}

\begin{theorem}[\rm cf. \cite{R1}]
\label{pi}
Let $M^n$ be a closed manifold with $\pi_1(M) =\pi$ and $\pi_i(M)=0$ for $2\le 
i 
\le n-2$. Then $\pi_{n-1}(M)=H^1_c(\pi)$. 
\end{theorem}

\p Let $\wt M$ be the universal covering space for $M$. Because of the 
Poincar\'e duality we have
$$
H^1_c(\pi)=H^1_c(\wt M)=H_{n-1}(\wt M).
$$
But, by the Hurewicz Theorem, $H_{n-1}(\wt M) = \pi_{n-1}(\wt 
M)=\pi_{n-1}(M)$.
\hqed

\begin{lemma}\label{null}
If $\pi$ is finitely generated nilpotent group with $\rk \pi
>1$, then $H^1_c(\pi)=0$.
\end{lemma}We embed $\pi$ as a 
uniform lattice in a contractible nilpotent Lie group $G$ with $\dim G = \rk 
\pi=n$, \cite{M}.

\p First, we assume that $\pi$ is torsion free. We embed $\pi$ as a 
uniform lattice in a contractible nilpotent Lie group $G$ with $\dim G = \rk 
\pi=n$, \cite{M}. Since $n>1$, we conclude that $\pi_{n-1}(G)=0$. Thus, 
since $\pi_1(G/\pi)=\pi$, we deduce from Theorem \ref{pi} that 
$H^1_c(\pi)=0$.

\m Now, if $\pi$ is not torsion free then it contains a torsion free subgroup 
$\pi'$ of finite index, \cite{Ku}. Then $K(\pi',1)$ can be regarded as a 
finite covering over $K(\pi,1)$. So, $K(\pi',1)$ and  $K(\pi,1)$ have 
the same universal covering, and thus, $H^1_c(\pi)=H^1_c(\pi')=0$. 
\hqed

\begin{cory}\label{rank}
Let $M$ be a closed $n$-dimensional manifold, $n>1$ with $\pi_i(M)=0$ for 
$i=2, 
\ldots, n-2$. If $\pi_1(M)$ is a nilpotent group $\pi$ with $\rk \pi >1$, then 
$\pi$ is torsion free and $\rk \pi =n$. 
\end{cory}

\p By Lemma \ref{null}, $H^1_c(\pi)=0$. Therefore, by \ref{pi}, 
$\pi_{n-1}(M)=0$. Furthermore, $H_n(\wt M)=0$ because $\pi$ is infinite. So, 
$\pi_n(\wt M)=0$, i.e. $\wt M$ is contractible, i.e.  $M=K(\pi,1)$. So, 
$\pi$ is torsion free since $M$ is finite dimensional. Finally, $M$ is 
homotopy equivalent to a closed nilmanifold $G/\pi$ of dimension $n$, and 
therefore $\rk \pi =n$.
\hqed 

\begin{theorem}
\label{4-nil}
Let $M$ be a closed $4$-dimensional symplectic manifold $M$ with $\pi_2(M)=0$. 
If $\pi_1(M)$ is a nilpotent, then $\pi_1(M)$ is a torsion free nilpotent 
group 
of rank $4$ and $H^4(\pi;\RR)=0$. Conversely, every finitely presented torsion 
free nilpotent group $\pi$ with $\rk \pi=4$ and $H^4(\pi)=\ZZ$ can be 
realized as the fundamental group of closed $4$-dimensional symplectically 
aspherical manifold $M$ with $\pi_2(M)=0$.
\end{theorem}

\p First, notice that $\rk \pi_1(M) >1$. Indeed, if $\rk \pi_1(M) =1$ then 
$\pi_1(M)$ contains $\ZZ$ as a subgroup of finite index. Considering the 
finite 
covering with respect to the inclusion $\ZZ \subset \pi_1(M)$, we get a 
4-dimensional closed symplectic manifold with the fundamental group $\ZZ$. But 
this is impossible by Proposition \ref{dim}.

\m Now, by Lemma \ref{rank}, $\pi_1(M)$  must be a torsion free nilpotent 
group 
of the rank 4. Furthermore, $H^4(\pi)=H^4(M)=\ZZ$, since $M$ is symplectic and 
therefore orientable.

\m Conversely, consider a torsion free finitely presented nilpotent group 
$\pi$, 
$\rk \pi=4$ an $H^4(\pi)=\ZZ$. We embed$\pi$ as a 
uniform lattice in a contractible nilpotent Lie group $G$ with $\dim G = \rk 
\pi=4$, \cite{M}. Then $M=G/pi$ is a 4-dimensional nilmanifold with 
$H^4(M)=\ZZ$. Now the result follows from Corollary~\ref{4nil-aspher}.
\hqed

\begin{cory}
Let $\pi$ be a torsion free finitely generated c-symplectic nilpotent group. 
If 
$\rk \pi > 4$ then $\pi\in \AAA_6, \pi\in \BBB_4$ and $\pi \notin \AAA_4$.
\end{cory}

This Corollary strength Corollary \ref{six}.

\p Recall that $\pi$ is a uniform lattice in a certain simply connected 
nilpotent Lie group 
$G$, $\dim G =\rk \pi$. Asserting as in Lemma \ref{c-sympl}, we conclude that 
$\pi\in 
\AAA_{2n}$. Therefore $\pi \in \AAA_6$, see Proposition \ref{chain}. 
Furthermore, $\pi\notin \AAA_4$ by Theorem \ref{4-nil}. Thus, $\pi\in \BBB_4$ 
since, by Proposition \ref{chain}, $\AAA_6\subset \AAA_4\cup \BBB_4$.
\hqed


\begin{thebibliography}{[IRTU]}\bibliographystyle{amsalpha}
          
\bibitem[B]{B}
K. Brown, {\it Cohomology of groups},
Springer-Verlag, Berlin, 1982.

\bibitem[D]{D}
S. Donaldson, Symplectic submanifolds and almost-complex geometry,
{\it J. Diff. Geom.}
{\bf 44} (1996) 666--705.

\bibitem[FLS]{FLS}
M. Fern\'andez, M. de Le\'on and M. Saralegui, A six-dimensional compact 
symplectic solvmanifold without K\"ahler
structures, 
{\it Osaka J. Math.}
{\bf  33}
(1996) 19--35.

\bibitem[F]{F}
A. Floer, Symplectic fixed points and holomorphic spheres,
{\it Comm. Math. Phys.}
{\bf  120}
(1989)
575--611.

\bibitem[G1]{G1}
R. Gompf, A new construction of symplectic manifolds, 
{\it Annals of Math.}
{\bf  142}
(1995) 
527--595.

\bibitem[G2]{G2}
R. Gompf, On symplectically aspherical symplectic manifolds with nontrivial 
$\pi_2$,
{\it Math. Res. Letters}
{\bf  5}
(1999) 
599--603.

\bibitem[Gb]{Gb}
V. V. Gorbatsevich, On lattices in Lie groups of types $(E)$ and $(R)$, 
{\it Moscow Univ. Math. Bull}.
{\bf 30} (1975) 98--104.

\bibitem[Ha]{Ha}
A. Hattori, Spectral sequence in de Rham cohomology of fibre bundles,
{\it J. Fac. Sci. Univ. Tokyo}
{\bf  8}  (Sect.1) (1960)  289--331.


\bibitem[H]{H}
H. Hofer, Lusternik--Schnirelmann theory for Lagrangian intersectons,
{\it Ann. Inst. H. Poincar\'e - Anal. Non Lineare}
{\bf 5} (1988)  465--499.

\bibitem[He]{He}
S. Helgason,
{\it Differential Geometry and Symmetric Spaces},
Academic Press, New York 1962.

\bibitem[IRTU]{IRTU}
R. Ib\'a\~nez, Yu. Rudyak, A. Tralle and L. Ugarte, On symplectically harmonic 
forms on six-dimensional nilmanifolds,
{\it Comment. Math. Helv.} {\bf  76}, no 1 (2001) 89--109, Erratum {\bf 76} 
(2001) 576.

\bibitem[K]{K}
J. K\c edra, Remarks on the flux group,
{\it Math. Res. Letters}
{\bf  7} (2000) 279--285.

\bibitem[Ku]{Ku}
A. G. Kurosh, {\it The theory of groups}, Chelsea
Publishing Co., New York 1960

\bibitem[LMP]{LMP}
F. Lalonde, D. McDuff\ and\ L. Polterovich, On the flux conjectures, 
in
{\it Geometry, topology, and dynamics {\rm (}Montreal, PQ, 1995{\rm )}}, 69--85, 
Amer. Math. Soc., Providence, RI 1998

\bibitem[LO]{LO}
G. Lupton, J. Oprea, Cohomologically symplectic spaces: toral actions and the 
Gottlieb group,
{\it Trans. Amer. Math. Soc.}
{\bf  347}, 1 (1995) 261--288.

\bibitem[M]{M}
A. I. Mal'cev, On a class of homogeneous spaces, {\it Izv AN SSSR 
Ser.  Matem.} {\bf 3} (1949) 9--32.
\bibitem[MS]{MS}
D. McDuff and D. Salamon
{\it Introduction to Symplectic Topology}, 
Clarendon Press, Oxford 1998.

\bibitem[N]{N}
S. P. Novikov, On manifolds with free abelian fundamental group 
and applications 
(Pontrjagin classes, smoothing, high--dimensional knots), {\it Izvestiya AN 
SSSR, 
Ser. Math.} 30 (1966) 208--246.

\bibitem[Ono]{Ono} 

K. Ono, Obstruction to circle group actions preserving symplectic structure, 
{\it Hokkaido Math. J.} {\bf 21} (1992), 99-102.

\bibitem[R1]{R1}Yu. B. Rudyak,
 On the fundamental group of a three-dimensional manifold, {\it Soviet Math. 
Doklady} {\bf 14} (1973)  814--818.
\bibitem[R2]{R2}
Yu. B. Rudyak, On strict category weight, gradient-like flows and the Arnold 
conjecture, 
{\it Internat. Math. Research Notices}
{\bf 5} (2000) 271--279.
 
\bibitem[RO]{RO}
Yu. Rudyak, J. Oprea, On the Lusternik--Schnirelmann category of symplectic 
manifolds and the Arnold conjecture,
{\it Math. Z.} 
{\bf  230}, no 4 (1999) 673--678.

\bibitem[RT]{RT}
Y. Rudyak and A. Tralle, On symplectic manifolds with aspherical symplectic 
form, 
{\it Topol. Methods Nonlinear Anal.}
{\bf  14} (1999) 353--362.

\bibitem[Sa]{Sa}
S. Salamon, Complex structures on nilpotent Lie algebras, {\it J. Pure Appl. 
Algebra} {\bf 157} (2001), no. 2-3, 311--333. 

\bibitem[St]{St}
N. Steenrod, {\it The topology of fibre bundles}, Reprint of the 1957 edition. 
Princeton Landmarks in Mathematics. Princeton
  Paperbacks. Princeton University Press, Princeton, NJ, 1999
  
\bibitem[TO]{TO}
A. Tralle and J. Oprea, {\it Symplectic manifolds with no K\"ahler structure}, 
Lect. Notes Math. 1661, Springer, Berlin 1997.

\bibitem[VGS]{VGS}
E.B. Vinberg, V.V. Gorbatsevich and O.V. Schwartsman, {\it Discrete subgroups 
in Lie groups}, 
in: A.L. Onishchik, E.B. Vinberg (eds), {Lie groups and Lie Algebras II}, 
1--123, 217--223, 
Encycl. Math. Sci. 21, Springer, Berlin 2000.

\end{thebibliography}
\enddocument
\end